\documentclass{article}

\usepackage{amssymb}
\usepackage{amsmath}

\usepackage[dvips]{graphicx}

\begin{document}


\begin{center}
\textbf{ON UNIFORM TOPOLOGICAL ALGEBRAS }
\end{center}
 
\noindent \textbf{}

\begin{center}
M. EL AZHARI 
\end{center}

\noindent \textbf{}

\noindent \textbf{Abstract.} The uniform norm on a uniform normed Q-algebra is the only
uniform Q-algebra norm on it. The uniform norm on a regular uniform normed
Q-algebra with unit is the only uniform norm on it. Let A be a uniform
topological algebra whose spectrum M(A) is equicontinuous, then A is a uniform
normed algebra. Let A be a regular semisimple commutative Banach algebra, then
every algebra norm on A is a Q-algebra norm on A.
 
\noindent \textbf{}


\noindent \textbf{} 


\noindent \textbf{}
 
\noindent \textbf{1. Preliminaries }                                                                        

\noindent \textbf{}

A topological algebra is an algebra (over the
complex field) which is also a Hausdorff topological vector space such that the
multiplication is separately continuous. A locally convex algebra is a
topological algebra whose topology is locally convex. A uniform seminorm on an
algebra $A$ is a seminorm $p$ satisfying $(i):p(x^{2})= p(x)^{2}$   
for all $x\in A.$ It is shown in [3], that a
seminorm satisfying $(i)$ is submultiplicative. A uniform topological algebra
(uT-algebra) is a topological algebra whose topology is determined by a family of
uniform seminorms. A uniform normed algebra is a normed algebra $(A,\Vert .\Vert)$ such
that  $\Vert x^{2}\Vert = \Vert x\Vert^{2}$ for all $x\in A.$ A topological algebra is a Q-algebra
if the set of quasi-invertible elements is open. A norm $\Vert .\Vert$ on an algebra 
$A$ is said to be an algebra norm (resp. Q-algebra norm) if $(A,\Vert .\Vert)$ is a normed
algebra (resp. normed Q-algebra). Let $A$ be an algebra, if $x\in A$ we denote
by $sp_{A}(x)$ the spectrum of $x$ and by $r_{A}(x)$ the spectral radius of $x.$  
An algebra $A$ is spectrally bounded if for each $x\in A, sp_{A}(x)$ is
bounded. For a topological algebra A; $M(A)$ denotes the set of all nonzero
continuous multiplicative linear functionals on $A, M(A)$ is endowed
with the weak topology induced by the topological dual of $A.$ For an arbitrary
topological algebra, $M(A)$ may be empty. A topological algebra is regular if
given a closed subset $F\subset M(A)$ and $f_{1}\in M(A), f_{1}\notin F,$   
there exists an $x\in A$ such that $\widehat{x}_{ \vert F}= 0$ and 
$\widehat{x}(f_{1})\neq 0; \widehat{x}:M(A)\rightarrow C, \widehat{x}(g)= g(x),$
is the Gelfand transform of $x.$ Let
$(B,\Vert .\Vert)$ be a uniform Banach algebra,  
$r_{B}(x)= \lim _{n\rightarrow\infty}\Vert x^{n}\Vert^{\frac{1}{n}}=\lim_{n\rightarrow\infty}\Vert x^{2^{n}} \Vert^{\frac{1}{2^{n}}} = \Vert x\Vert $
for  all $x\in B; B$ is commutative by [4, Corollary 1], hence
$\Vert x\Vert = sup\lbrace\vert f(x)\vert, f\in M(B)\rbrace $ for all $x\in B.$ Let $A$ be a
uT-algebra; the completion of $A$ is a uT-algebra and then an inverse limit of
uniform Banach algebras, thus $A$ is commutative.
 
\noindent \textbf{}

\noindent \textbf{2. Results}

\noindent \textbf{}

\noindent \textbf{Theorem 2.1.} Let $(A,\Vert .\Vert)$ be a uniform normed Q-algebra. If $\vert .\vert$ is
a uniform norm on $A$, then $\vert x\vert\leq\Vert x\Vert$ for all $x\in A.$

\noindent \textbf{}
 
\noindent \textbf {Proof.} By [7, Proposition 7.5], $M(A,\Vert .\Vert)$(resp. $
M(A,\vert .\vert)$) is topologically isomorphic to $M(\tilde{A},\Vert .\Vert)$(resp.
$M(\tilde{A},\vert .\vert)),(\tilde{A},\Vert .\Vert)$ and $(\tilde{A},\vert .\vert)$ are
respectively the completions of $(A,\Vert .\Vert)$ and $(A,\vert .\vert). (\tilde{A},\Vert .\Vert)$ and $(\tilde{A},\vert .\vert)$ are uniform Banach algebras,
then for $x\in A,$ we have $\Vert x\Vert =\sup\lbrace\vert f(x)\vert, 
f\in M(\tilde{A},\Vert .\Vert)\rbrace = \sup\lbrace\vert f(x)\vert, f\in M(A,\Vert .\Vert)\rbrace$
and $\vert x\vert =\sup\lbrace\vert f(x)\vert, 
f\in M(\tilde{A},\vert .\vert)\rbrace = \sup\lbrace\vert f(x)\vert, f\in M(A,\vert .\vert)\rbrace$,                       
hence $ \vert x\vert\leq\Vert x\Vert$ since
$M(A,\vert .\vert)\subset M(A,\Vert .\Vert)$; the last inclusion is due to the fact
that $(A,\Vert .\Vert)$ is a Q-algebra.

\noindent \textbf{}

\noindent \textbf{Corollary 2.2.} Let $(A,\Vert .\Vert)$ be a uniform normed
Q-algebra, then $\Vert .\Vert$ is the unique uniform Q-algebra norm on $A.$

\noindent \textbf{}
 
\noindent \textbf{Theorem 2.3.} Let $(A,\Vert .\Vert)$ be a regular uniform normed
Q-algebra with unit. If $\vert .\vert$ is a uniform norm on $A,$ then $ \Vert .\Vert =\vert .\vert.$

\noindent \textbf{}

\noindent \textbf{Proof.} Let $ x\in A,\vert x\vert =\sup\lbrace\vert f(x)\vert, 
f\in M(\tilde{A},\vert .\vert)\rbrace = \sup\lbrace\vert f(x)\vert, f\in M(A,\vert .\vert)\rbrace$.     
$M(A,\vert .\vert)$ is topologically isomorphic to $M(\tilde{A},\vert .\vert)$,
then  $M(A,\vert .\vert)$ is compact. Put $T=M(A,\vert .\vert),T$ is a closed subset of $ 
M(A,\Vert .\Vert)$. Suppose that $T\neq M(A,\Vert .\Vert)$, let 
$f_{1}\in M(A,\Vert .\Vert)$ such that $f_{1}\notin T$. Since
$(A,\Vert .\Vert)$ is regular, there exists  $x\in A$ such that $\widehat{x}_{ \vert T}=0$ and 
$\widehat{x}(f_{1})\neq 0$, hence $ \vert x\vert= 0$ with  
$x\neq 0$, a contradiction.  For $x\in A, \Vert x\Vert =\sup\lbrace\vert f(x)\vert, 
f\in M(\tilde{A},\Vert .\Vert)\rbrace = \sup\lbrace\vert f(x)\vert, f\in M(A,\Vert .\Vert)\rbrace = \sup\lbrace\vert f(x)\vert, f\in M(A,\vert .\vert)\rbrace =\vert x\vert.$   
 
\noindent \textbf{}
 
\noindent \textbf{Theorem 2.4}. Let $A$ be a uT-algebra whose spectrum $M(A)$ is
equicontinuous. Then $A$ is a uniform normed algebra.

\noindent \textbf{}
 
\noindent \textbf{Proof.} The topology of $A$ is determined by a family $\lbrace p_{u}, 
u\in U\rbrace$ of uniform seminorms. For each $u\in U$, let 
$ N_{u}=\lbrace x\in A,p_{u}(x)=0\rbrace$ and let $A_{u}$ be the Banach algebra obtained by
completing $A/N_{u}$ in the norm $ \Vert x_{u}\Vert_{u}= p_{u}(x),x_{u}= x+N_{u}$, 
it is clear that $A_{u}$ is a uniform Banach algebra. For each 
$u\in U$, let $M_{u}(A)= \lbrace f\in M(A),\vert f(x)\vert\leq p_{u}(x)$     
for all $x\in A\rbrace$.  
Let $u\in U$ and $x\in A, p_{u}(x)=\Vert x_{u}\Vert_{u}=\sup\lbrace\vert g(x_{u})\vert,g\in M(A_{u})\rbrace=\sup\lbrace\vert f(x)\vert,f\in M_{u}(A)\rbrace$ by [7,Proposition 7.5].     
Let $q(x)=\sup\lbrace\vert f(x)\vert,f\in M(A)\rbrace,x\in A.$ Since $M(A)$ is equicontinuous, it follows
that $q$ is a continuous seminorm on $A$. Let $x\in A$ and suppose that
$q(x)=0,$ since $A$ is Hausdorff and 
$p_{u}(x)\leq q(x)=0$ for all $u\in U,$ it follows that
$x=0.$ q is a continuous uniform norm on $A,$ also
$p_{u}(x)\leq q(x)$ for all $u\in U$ and $x\in A,$ then
the topology of $A$ can be defined by the uniform norm $q.$
 
\noindent \textbf{}

\noindent \textbf{Corollary 2.5.} Let $A$ be a uT-algebra that is a Q-algebra. Then $A$ is
a uniform normed algebra.
 
\noindent \textbf{}

\noindent \textbf{Proof.} By [6, Proposition II.7.1], every topological Q-algebra has
an equicontinuous spectrum.
 
\noindent \textbf{}

\noindent \textbf{Corollary 2.6.} Let $A$ be a spectrally bounded, barrelled, uT-algebra. 
Then $A$ is a uniform normed algebra.
 
\noindent \textbf{}

\noindent \textbf{Proof.} $\sup\lbrace\vert f(x)\vert,f\in M(A)\rbrace\leq r_{A}(x)$ 
for all $x\in A. M(A)$ is bounded for the weak topology, then
$M(A)$ is equicontinuous.
 
\noindent \textbf{}

\noindent \textbf{Remark.} These results show that completeness is not necessary for some
results of S. J. Bhatt on uniform topological algebras (see  [1, Theorem 1(ii)],[1, the second affirmation in the corollary],[1, Theorem 2] and [2, Corollary 2.5].
 
\noindent \textbf{}

\noindent \textbf{Theorem 2.7.} Let $A$ be a topological algebra. If $p$ is a uniform
seminorm on $A$, then $p(x)\leq r_{A}(x)$ for all $x\in A.$ If $A$ is a Q-algebra, 
then every uniform seminorm on $A$ is continuous.
 
\noindent \textbf{} 

\noindent \textbf{Proof.} Let $N_{p}=\lbrace x\in A,p(x)=0\rbrace$ and $A_{p}$ be the
Banach algebra obtained by completing $A/N_{p}$ in the norm 
$\Vert x_{p}\Vert_{p}=p(x), x_{p}=x+N_{p}. \Vert .\Vert_{p}$ is a uniform norm on $A_{p}$
since $p$ is a uniform seminorm on $A.$ Consider $G:A\rightarrow A_{p},G(x)=x_{p}$,
the quotient map.  Let $x\in A, p(x)=\Vert x_{p}\Vert_{p}=r_{A_{p}}(x_{p})\leq r_{A}(x).$    
If $A$ is a Q-algebra, $r_{A}$ is continuous at 0 by [6, Lemma II.4.2], then $p$ is
continuous.
 
\noindent \textbf{}

\noindent \textbf{Remark.} Theorem 2.7 shows that the hypothesis  "$A$ is a unital locally
convex algebra and the inversion in $A$ is continuous" is not necessary in [2, Theorem 2.3].
 
\noindent \textbf{} 

Now we generalize the first affirmation of [1, Corollary].
 
\noindent \textbf{} 

\noindent \textbf{Theorem 2.8.} Let $(A,\Vert .\Vert)$ be a regular semisimple
commutative Banach algebra. Then every algebra norm on $A$ is a Q-algebra norm on
$A.$
 
\noindent \textbf{} 

\noindent \textbf{Proof.} (1) $A$ is unital: Let $\vert .\vert$ be an algebra norm on $A.$ Put 
$K_{1}=M(A,\vert .\vert)$ and $K=M(A,\Vert .\Vert), K_{1}\subset K, K_{1}$ is
homeomorphic to $M(\tilde{A},\vert .\vert) ((\tilde{A},\vert .\vert)$
is the completion of $(A,\vert .\vert)),$ then $K_{1}$ is compact, hence 
$K_{1}$ is closed in $K.$ Suppose that $K_{1}\neq K. K\backslash K_{1}$ is
a nonempty open set, there exists an open set $G$ in $K$ such that 
$\bar{G}\subset K\backslash K_{1}$ since $(A,\Vert .\Vert)$ is regular. 
Also by regularity of $(A,\Vert .\Vert)$ [5, Corollary 7.3.4], there exist
$x\in A,y\in A,y\neq 0,$ such that $\widehat{x}_{\vert K_{1}}=1,\widehat{x}_{\vert\bar{G}}=0$   
and $\widehat{y}_{\vert K\backslash G}=0$.
$x$ is invertible in $(\tilde{A},\vert .\vert)$ since $\widehat{x}_{\vert K_{1}}=1$. 
$\widehat{xy}=\widehat{x}\widehat{y}=0$ on $K$, then $xy\in Rad(A)=\lbrace 0\rbrace$, 
hence $y=x^{-1}xy=0,$ a contradiction.  Finally, $K_{1}=K.$ Let $x\in A,r_{A}(x)=\sup \lbrace\vert f(x)\vert,f\in K\rbrace=\sup \lbrace\vert f(x)\vert,f\in K_{1}\rbrace\leq \vert x\vert$,  
then $(A,\vert .\vert)$ is a Q-algebra by [8, Lemma 2.1].    
\newline (2) $A$ is not unital:  Let $A_{1}$ be the algebra obtained from $A$ by adjunction of an identity $e.$  The elements of  $A_{1}$ have
the form  $x+re,$ where $\ x\in A$ and $r\in C$. $A_{1}$ is also regular and
semisimple. Let $p$ be an algebra norm on $A.$ For $x\in A$ and $r\in C,$ put 
$p_{1}(x+re)=p(x)+\vert r\vert,$ it is easy to show that $p_{1}$
is an algebra norm on $A_{1}.$ By (1) and [8, Lemma 2.1],  $r_{A_{1}}(x+re)\leq 
p_{1}(x+re)$ for all $x\in A$ and $r\in C,$ consequently  
$r_{A}(x)=r_{A_{1}}(x)\leq p_{1}(x)=p(x)$ for
all $x\in A.$ Then $p$ is a Q-algebra norm on $A$ by [8, Lemma 2.1].
 
\noindent \textbf{}  

The results of this paper were communicated at the international conference on topological algebras
and applications (ICTAA 2000).

\noindent \textbf{}

\begin{center}
REFERENCES
\end{center}
 
\noindent \textbf{} 
 
\noindent \textbf{} [1] S. J. Bhatt and D. J. Karia, Uniqueness of the uniform norm with application
to topological algebras, Proc. Amer. Math. Soc., 116 (1992), 499-503.

\noindent \textbf{} [2] S. J. Bhatt, Automatic continuity of homomorphisms in topological algebras,
Proc. Amer. Math. Soc., 119 (1993), 135-139.

\noindent \textbf{} [3] H. V. Dedania, A seminorm with square property is automatically
submultiplicative, Proc. Indian. Acad. Sci. (Math. Sci.), 108 (1998), 51-53.

\noindent \textbf{} [4] R. A. Hirschfeld and W. Zelazko, On spectral norm Banach algebras, Bull.
Acad. Pol. Sc., XVI  (1968), 195-199.

\noindent \textbf{} [5] R. Larsen, Banach algebras, Marcel-Dekker, New-York, 1973.

\noindent \textbf{} [6] A. Mallios, Topological algebras, Selected Topics, North-Holland, Amsterdam,
1986.

\noindent \textbf{} [7] E. A. Michael, Locally multiplicatively convex topological algebras, Mem.
Amer. Math. Soc.,  11 (1952).

\noindent \textbf{} [8] B. Yood, Homomorphisms on normed algebras, Pacific J. Math.,  8 (1958),
373-381.

\noindent \textbf{} 

\noindent \textbf{} Ecole Normale Sup\'{e}rieure

\noindent \textbf{} Avenue Oued Akreuch

\noindent \textbf{} Takaddoum, BP 5118, Rabat

\noindent \textbf{} Morocco
 
\noindent \textbf{} 

\noindent \textbf{} E-mail:  mohammed.elazhari@yahoo.fr

\end{document}